\documentclass[12pt]{article}

\usepackage{amssymb,amsbsy,amsthm,amsmath,amscd,epsfig,graphicx,times}

\input{arrow}

\begin{document}

\title{\Large{\textbf{DEDUCING THE MULTIDIMENSIONAL SZEMER\'EDI THEOREM FROM AN INFINITARY REMOVAL LEMMA}}}
\author{Tim Austin}
\date{}

\maketitle


\newenvironment{nmath}{\begin{center}\begin{math}}{\end{math}\end{center}}

\newtheorem{thm}{THEOREM}[section]
\newtheorem*{thm*}{THEOREM}
\newtheorem{lem}[thm]{LEMMA}
\newtheorem{prop}[thm]{PROPOSITION}
\newtheorem{cor}[thm]{COROLLARY}
\newtheorem*{conj*}{CONJECTURE}
\newtheorem{dfn}[thm]{DEFINITION}
\newtheorem{ques}[thm]{QUESTION}
\theoremstyle{remark}


\newcommand{\A}{\mathcal{A}}
\newcommand{\B}{\mathcal{B}}
\newcommand{\I}{\mathcal{I}}
\newcommand{\frH}{\mathfrak{H}}
\renewcommand{\Pr}{\mathrm{Pr}}
\newcommand{\s}{\sigma}
\renewcommand{\P}{\mathcal{P}}
\renewcommand{\O}{\Omega}
\renewcommand{\S}{\Sigma}
\newcommand{\T}{\mathrm{T}}
\newcommand{\co}{\mathrm{co}}
\newcommand{\e}{\mathrm{e}}
\newcommand{\eps}{\varepsilon}
\renewcommand{\d}{\mathrm{d}}
\newcommand{\im}{\mathrm{i}}
\renewcommand{\l}{\lambda}
\newcommand{\U}{\mathcal{U}}
\newcommand{\G}{\Gamma}
\newcommand{\g}{\gamma}
\newcommand{\calL}{\mathcal{L}}
\renewcommand{\L}{\Lambda}
\newcommand{\hcf}{\mathrm{hcf}}
\newcommand{\FLat}{\mathrm{FLat}}
\newcommand{\F}{\mathcal{F}}
\renewcommand{\a}{\alpha}
\newcommand{\bbN}{\mathbb{N}}
\newcommand{\bbR}{\mathbb{R}}
\newcommand{\bbZ}{\mathbb{Z}}
\newcommand{\bbQ}{\mathbb{Q}}
\newcommand{\bbT}{\mathbb{T}}
\newcommand{\sfE}{\mathsf{E}}
\newcommand{\sfP}{\mathsf{P}}
\newcommand{\id}{\mathrm{id}}
\newcommand{\bb}[1]{\mathbb{#1}}
\newcommand{\fr}[1]{\mathfrak{#1}}
\renewcommand{\bf}[1]{\mathbf{#1}}
\renewcommand{\rm}[1]{\mathrm{#1}}
\renewcommand{\cal}[1]{\mathcal{#1}}
\newcommand{\fin}{\nolinebreak\hspace{\stretch{1}}$\lhd$}
\newcommand{\uhr}{\!\!\upharpoonright}

\begin{abstract}
We offer a new proof of the Furstenberg-Katznelson multiple
recurrence theorem for several commuting probability-preserving
transformations $T_1$, $T_2$, \ldots, $T_d:\bbZ\curvearrowright
(X,\S,\mu)$~(\cite{FurKat78}), and so, via the Furstenberg
correspondence principle introduced in~\cite{Fur77}, a new proof of
the multi-dimensional Szemer\'edi Theorem.  We bypass the careful
manipulation of certain towers of factors of a
probability-preserving system that underlies the
Furstenberg-Katznelson analysis, instead modifying an approach
recently developed in~\cite{Aus--nonconv} to pass to a large
extension of our original system in which this analysis greatly
simplifies. The proof is then completed using an adaptation of
arguments developed by Tao in~\cite{Tao07} for his study of an
infinitary analog of the hypergraph removal lemma. In a sense, this
addresses the difficulty, highlighted by Tao, of establishing a
direct connection between his infinitary, probabilistic approach to
the hypergraph removal lemma and the infinitary, ergodic-theoretic
approach to Szemer\'edi's Theorem set in motion by
Furstenberg~\cite{Fur77}.
\end{abstract}

\parskip 0pt

\tableofcontents

\parskip 7pt

\section{Introduction}

We give a new ergodic-theoretic proof of the multidimensional
multiple recurrence theorem of Furstenberg and
Katznelson~\cite{FurKat78}, which their correspondence principle
shows to be equivalent to the multidimensional Szemer\'edi Theorem.

\begin{thm}[Multidimensional multiple recurrence]\label{thm:multirec}
Suppose that $T_1,T_2,\ldots,T_d:\bbZ\curvearrowright (X,\S,\mu)$
are commuting probability-preserving actions and that $A \in \S$ has
$\mu(A) > 0$. Then
\[\liminf_{N\to\infty}\frac{1}{N}\sum_{n=1}^N\mu(T_1^{-n}(A)\cap T_2^{-n}(A)\cap\cdots\cap T_d^{-n}(A)) > 0,\]
and so, in particular, there is some $n \geq 1$ with
\[\mu(T_1^{-n}(A)\cap T_2^{-n}(A)\cap\cdots\cap T_d^{-n}(A)) > 0.\]
\end{thm}

Our proof of Theorem~\ref{thm:multirec} will call on some rather
different ergodic-theoretic machinery from Furstenberg and
Katznelson's.  Our main technical ingredients are the notions of
`pleasant' and `isotropized' extensions of a system. Pleasant
extensions were first used in~\cite{Aus--nonconv} to give a new
proof of the (rather easier) result that the `nonconventional
ergodic averages'
\begin{eqnarray}\label{eq:nonconv}
\frac{1}{N}\sum_{n=1}^N\prod_{i\leq d}f_i\circ T_i^n
\end{eqnarray}
associated to $f_1,f_2,\ldots,f_d \in L^\infty(\mu)$ always converge
in $L^2(\mu)$ as $N\to\infty$.  (This was first shown by Tao
in~\cite{Tao08(nonconv)}, although various special cases had
previously been established by other
methods~\cite{ConLes88.1,ConLes88.2,Zha96,HosKra01,HosKra05,Zie07,FraKra05}.)
Much of the present paper is motivated by the results used
in~\cite{Aus--nonconv} to give a new proof of this convergence.
Isotropized extensions are a new tool developed for the present
paper, but their analysis is closely analogous to that of pleasant
extensions.

After passing to a pleasant and isotropized extension of our
original system, the limit of~(\ref{eq:nonconv}) takes a special
form, and in this paper it is by analyzing this expression that we
shall prove positivity. It turns out that this special form enables
us to make contact with the machinery developed by Tao
in~\cite{Tao07} for his infinitary proof of the hypergraph removal
lemma.  Since the hypergraph removal lemma offers a known route to
proving the multidimensional Szemer\'edi Theorem (as shown, subject
to some important technical differences, by Nagle, R\"odl and
Schacht~\cite{NagRodSch07} and by Gowers~\cite{Gow??}), and this in
turn is equivalent to multidimensional multiple recurrence, Tao's
work already offers a proof of multiple recurrence using his
infinitary removal lemma. In a sense, our present contribution is to
short-circuit the above chain of implications and give a near-direct
proof of multiple recurrence using Tao's ideas.  Unfortunately, we
have not been able to make a reduction to a simple black-box appeal
to Tao's result; rather, we formulate
(Proposition~\ref{prop:infremoval}) a closely-related result adapted
to our ergodic theoretic setting, which then admits a very similar
proof. With this caveat, our work addresses the question of relating
infinitary proofs of multiple recurrence and hypergraph removal
explicitly raised by Tao at the beginning of Section 5
of~\cite{Tao07}: it turns out that his ideas are not directly
applicable to an arbitrary probability-preserving $\bbZ^d$-system,
but becomes so only when we enlarge the system to lie in the special
class of systems that are pleasant and isotropized.

\textbf{Acknowledgements}\quad My thanks go to Vitaly Bergelson and
David Fremlin for helpful comments on an earlier version of this
paper, and to the Mathematical Sciences Research Institute where a
significant re-write was undertaken after a serious flaw was
discovered in an earlier version.

\section{Basic notation and preliminaries}\label{sec:prelim}

Throughout this paper $(X,\S)$ will denote a measurable space.
Since our main results pertain only to the joint distribution of
countably many bounded real-valued functions on this space and their
shifts under some measurable transformation, by passing to the image
measure on a suitable shift space we may always assume that $(X,\S)$
is standard Borel, and this will prove convenient for some of our
later constructions. In addition, $\mu$ will always denote a
probability measure on $\S$. We shall write $(X^e,\S^{\otimes e})$
for the usual product measurable structure indexed by a set $e$, and
$\mu^{\otimes e}$ for the product measure and $\mu^{\Delta e}$ for
the diagonal measure on this structure respectively.  We also write
$\pi_i:X^e \to X$ for the $i^{\rm{th}}$ coordinate projection
whenever $i \in e$. Given a measurable map $\phi:(X,\S)\to (Y,\Phi)$
to another measurable space, we shall write $\phi\circ\mu$ for the
resulting pushforward probability measure on $(Y,\Phi)$.

If $T:\G\curvearrowright (X,\S,\mu)$ is a probability-preserving
action of a countable group $\G$, then by a \textbf{factor} of the
quadruple $(X,\S,\mu,T)$ we understand a globally $T$-invariant
sub-$\s$-algebra $\Phi \leq \S$. The \textbf{isotropy factor} is the
sub-$\s$-algebra of those subsets $A \in \S$ such that
$\mu(A\triangle T^\g(A)) = 0$ for all $\g\in \G$, and we shall
denote it by $\S^T$.  If $T_1,T_2:\G\curvearrowright (X,\mu)$ are
two commuting actions of the same \emph{Abelian} group, then we can
define another action by $(T_1^{-1}T_2)^\g := T^{\g^{-1}}_1T_2^\g$,
and then we write $\S^{T_1 = T_2}$ for $\S^{T_1^{-1}T_2}$, and
similarly if we are given a larger number of actions of the same
group. The most important kind of morphism from one $\G$-system
$T:\G\curvearrowright (X,\S,\mu)$ to another $S:\G\curvearrowright
(Y,\Phi,\nu)$ is given by a measurable map $\phi:X \to Y$ such that
$\nu = \phi\circ\mu$ and $S\circ\phi = \phi\circ T$: we call such a
$\phi$ a \textbf{factor map}.  In this case we shall write
$\phi:(X,\S,\mu,T) \to (Y,\Phi,\nu,S)$.  To a factor map $\phi$ we
can associate the factor $\{\phi^{-1}(A):\ A\in\Phi\}$.

Our specific interest is in $d$-tuples of commuting $\bbZ$-actions
$T_i$, $i=1,2,\ldots,d$.  Clearly these can be interpreted as the
$\bbZ$-subactions of a single $\bbZ^d$-action corresponding to the
$d$ coordinate directions $\bbZ\cdot e_i \leq \bbZ^d$.

Given these actions, we shall make repeated reference to certain
factors assembled from the isotropy factors among the $T_i$. These
will be indexed by subsets of $[d] := \{1,2,\ldots,d\}$, or more
generally by subfamilies of the collection $\binom{[d]}{\geq 2}$ of
all subsets of $[d]$ of size at least $2$.  On the whole, these
indexing subfamilies will be \textbf{up-sets} in $\binom{[d]}{\geq
2}$: $\I\subseteq \binom{[d]}{\geq 2}$ such that $u \in \I$ and
$[d]\supseteq v\supseteq u$ imply $v \in \I$. For example, given $e
\subseteq [d]$ we write $\langle e\rangle := \{u \in\binom{[d]}{\geq
2}:\ u\supseteq e\}$ (note the non-standard feature of our notation
that $e \in \langle e\rangle$ if and only if $|e| \geq 2$): up-sets
of this form are \textbf{principal}.  We will abbreviate
$\langle\{i\}\rangle$ to $\langle i\rangle$. It will also be helpful
to define the \textbf{depth} of a non-empty up-set $\I$ to be
$\min\{|e|:\ e\in \I\}$.

The corresponding factors are obtained for $e =
\{i_1,i_2,\ldots,i_k\}\subseteq [d]$ with $k\geq 2$ by defining
$\Phi_e := \S^{T_{i_1} = T_{i_2} = \ldots = T_{i_k}}$, and given an
up-set $\cal{I} \subseteq \binom{[d]}{\geq 2}$ by defining
$\Phi_{\cal{I}} := \bigvee_{e\in\cal{I}}\Phi_e$.

From the ordering among the factors $\Phi_e$ it is clear that
$\Phi_{\cal{I}} = \Phi_{\cal{A}}$ whenever $\cal{A} \subseteq
\binom{[d]}{\geq 2}$ is a family that generates $\cal{I}$ as an
up-set, and in particular that $\Phi_e = \Phi_{\langle e\rangle}$.

An \textbf{inverse system} is a family of probability-preserving
systems $T^{(m)}:\G\curvearrowright(X^{(m)},\S^{(m)},\mu^{(m)})$
together with factor maps
\[\psi_m:(X^{(m+1)},\S^{(m+1)},\mu^{(m+1)},T^{(m+1)})\to
(X^{(m)},\S^{(m)},\mu^{(m)},T^{(m)});\] from this one can construct
the \textbf{inverse limit}
\[\lim_{m\leftarrow}\,(X^{(m)},\S^{(m)},\mu^{(m)},T^{(m)})\] as
described, for example, in Section 6.3 of Glasner~\cite{Gla03}.

Finally, the following distributional condition for families of
factors will play a central r\^ole through this paper.

\begin{dfn}[Relative independence for factor-tuples]
If $\S_i \geq \Xi_i$ are factors of $(X,\S,\mu)$ for each $i\leq d$,
then the tuple of factors $(\S_1,\S_2,\ldots,\S_d)$ is
\textbf{relatively independent} over the tuple
$(\Xi_1,\Xi_2,\ldots,\Xi_d)$ if whenever $f_i\in
L^\infty(\mu\uhr_{\S_i})$ for each $i\leq d$ we have
\[\int_X\prod_{i\leq d}f_i\,\d\mu = \int_X\prod_{i\leq d}\sfE_\mu(f_i\,|\,\Xi_i)\,\d\mu.\]
\end{dfn}

\section{The Furstenberg self-joining}

It turns out that a particular $d$-fold self-joining of $\mu$ both
controls the convergence of the nonconventional
averages~(\ref{eq:nonconv}) and then serves to express their
limiting value.

Given our commuting actions and any $e = \{i_1 < i_2 < \ldots <
i_k\}\subseteq [d]$, we define
\begin{multline*}
\mu^{\rm{F}}_e(A_1\times A_2\times \cdots\times A_k)\\
:= \lim_{N\to
\infty}\frac{1}{N}\sum_{n=1}^N\mu(T_{i_1}^{-n}(A_1)\cap
T_{i_2}^{-n}(A_2)\cap\cdots\cap T_{i_k}^{-n}(A_k))
\end{multline*}
for $A_1,A_2,\ldots,A_k \in \S$.  That these limits always exist
(and so this definition is possible) follows from the convergence of
the nonconventional averages~(\ref{eq:nonconv}), although approaches
to convergence that use this self-joining (as
in~\cite{Aus--nonconv}, or for various special cases in~\cite{Zha96}
and~\cite{Zie07}) actually handle both kinds of limits alternately
in a combined proof of their existence by induction on $k$.

Given the existence of the limits (\ref{eq:nonconv}) and the
assumption that $(X,\S)$ is standard Borel, it is easy to check that
$\mu^{\rm{F}}_e$ extends to a $k$-fold self-joining of $\mu$ on
$\S^{\otimes e}$. This is the \textbf{Furstenberg self-joining} of
$\mu$ associated to $T_{i_1}$, $T_{i_2}$, \ldots, $T_{i_k}$. It is
now clear from our definition that the assertion of
Theorem~\ref{thm:multirec} can be re-stated as being that if $\mu
(A)
> 0$ then also $\mu^{\rm{F}}_{[d]}(A^d) > 0$. It is in this form
that we shall prove it.

The following elementary properties of the Furstenberg self-joining
will be important later.

\begin{lem}\label{lem:Fberg-project}
If $e = \{i_1 < i_2 < \ldots < i_k\} \subseteq e' = \{j_1 < j_2 <
\ldots < j_l\}$ then
$\pi_{\{i_1,i_2,\ldots,i_k\}}\circ\mu^{\rm{F}}_{e'} =
\mu^{\rm{F}}_e$.
\end{lem}

\textbf{Proof}\quad This is immediate from the definition: if
$A_{i_j} \in \S$ for each $j \leq k$ then
\begin{multline*}
(\pi_{\{i_1,i_2,\ldots,i_k\}}\circ\mu^{\rm{F}}_{e'})(A_1\times A_2\times \cdots\times A_k)\\
:= \lim_{N\to
\infty}\frac{1}{N}\sum_{n=1}^N\mu(T_{j_1}^{-n}(B_1)\cap
T_{j_2}^{-n}(B_2)\cap\cdots\cap T_{j_l}^{-n}(B_l))
\end{multline*}
where $B_j := A_j$ if $j \in e$ and $B_j := X$ otherwise; but then
this last average simplifies summand-by-summand directly to
\begin{multline*}
\lim_{N\to \infty}\frac{1}{N}\sum_{n=1}^N\mu(T_{j_1}^{-n}(A_1)\cap
T_{i_2}^{-n}(A_2)\cap\cdots\cap T_{i_k}^{-n}(A_k))\\
=: \mu^{\rm{F}}_e(A_1\times A_2\times \cdots\times A_k),
\end{multline*}
as required. \qed

\begin{lem}\label{lem:diag}
For any $e \subseteq [d]$ the restriction
$\mu^{\rm{F}}_e\uhr_{\Phi_e^{\otimes e}}$ is just the diagonal
measure $(\mu\uhr_{\Phi_e})^{\Delta e}$.
\end{lem}

\textbf{Proof}\quad If $e = \{i_1 < i_2 < \ldots < i_k\}$ and $A_j
\in \Phi_e$ for each $j \leq k$ then by definition we have
\begin{eqnarray*}\mu^{\rm{F}}_e(A_1 \times A_2\times \cdots\times A_k) &=& \lim_{N\to \infty}\frac{1}{N}\sum_{n=1}^N\mu(T_{i_1}^{-n}(A_1)\cap T_{i_2}^{-n}(A_2)\cap\cdots\cap
T_{i_k}^{-n}(A_k))\\
&=& \lim_{N\to
\infty}\frac{1}{N}\sum_{n=1}^N\mu(T_{i_1}^{-n}(A_1\cap
A_2\cap\cdots\cap A_k))\\
&=& \mu(A_1\cap A_2\cap\cdots\cap A_k),
\end{eqnarray*}
as required. \qed

It follows from the last lemma that whenever $e \subseteq e'$ the
factors $\pi_i^{-1}(\Phi_e) \leq \S^{\otimes e'}$ for $i\in e$ are
all equal up to $\mu^{\rm{F}}_{e'}$-negligible sets.  It will prove
helpful later to have a dedicated notation for these factors.

\begin{dfn}[Oblique copies]
For each $e \subseteq [d]$ we refer to the common
$\mu^{\rm{F}}_{[d]}$-completion of the sub-$\s$-algebra
$\pi_i^{-1}(\Phi_e)$, $i \in e$, as the \textbf{oblique copy} of
$\Phi_e$, and denote it by $\Phi^{\rm{F}}_e$.  More generally we
shall refer to factors formed by repeatedly applying $\cap$ and
$\vee$ to such oblique copies as \textbf{oblique factors}.
\end{dfn}

It will be important to know that Furstenberg self-joinings behave
well under inverse limits.  The following is another immediate
consequence of the definition, and we omit the proof.

\begin{lem}\label{lem:Fberg-inv-lim}
If
\[\ldots\to(X^{(m+1)},\S^{(m+1)},\mu^{(m+1)},T^{(m+1)})\stackrel{\psi_m}{\longrightarrow}(X^{(m)},\S^{(m)},\mu^{(m)},T^{(m)})\to\ldots\]
is an inverse system with inverse limit
$(\tilde{X},\tilde{\S},\tilde{\mu},\tilde{T})$, then the Furstenberg
self-joinings $\big((X^{(m)})^d,(\S^{(m)})^{\otimes
d},(\mu^{(m)})^{\rm{F}}_{[d]},T^{\times d}\big)$ with factor maps
$\phi_m^{\times d}$ also form an inverse system with inverse limit
$\big(\tilde{X}^d,\tilde{\S}^{\otimes
d},\tilde{\mu}^{\rm{F}}_{[d]},\tilde{T}^{\times d}\big)$. \qed
\end{lem}

\section{Pleasant and isotropized extensions}

We now introduce the main technical definitions of this paper: that
of `pleasant systems', closely following~\cite{Aus--nonconv}, and
alongside them the related notion of `isotropized systems'. Recall
that to a commuting tuple of actions
$T_1,T_2,\ldots,T_d:\bbZ\curvearrowright (X,\S,\mu)$ we associate
the factors
\[\Phi_e := \S^{T_{i_1} = T_{i_2} = \ldots = T_{i_k}}\]
indexed by subsets $e = \{i_1,i_2,\ldots,i_k\}\subseteq [d]$.

\begin{dfn}[Pleasant system]\label{dfn:pleasant}
A system $(X,\S,\mu,T)$ is \textbf{$(e,i)$-pleasant} for some $i \in
e\in\binom{[d]}{\geq 2}$ if the $i^{\rm{th}}$ coordinate projection
$\pi_i$ is relatively independent from the other $\pi_j$, $j \in e$,
over the factor $\pi_i^{-1}\big(\bigvee_{j\in e
\setminus\{i\}}\Phi_{\{i,j\}}\big)$ under the Furstenberg
self-joining $\mu^{\rm{F}}_e$:
\[\int_{X^e}\prod_{j\in e}f_j\circ\pi_j\,\d\mu^{\rm{F}}_e =
\int_{X^e}\Big(\sfE_\mu\Big(f_i\,\Big|\,\bigvee_{j\in e
\setminus\{i\}}\Phi_{\{i,j\}}\Big)\circ\pi_i\Big)\cdot\prod_{j\in
e\setminus \{i\}}f_j\circ\pi_j\,\d\mu^{\rm{F}}_e\] whenever $f_j \in
L^\infty(\mu)$ for each $j\in e$.

It is \textbf{fully pleasant} if it is $(e,i)$-pleasant for every
pair $i\in e$.
\end{dfn}

\begin{dfn}[Isotropized system]\label{dfn:isotropized}
A commuting tuple of actions
$T_1,T_2,\ldots,T_d:\bbZ\curvearrowright (X,\S,\mu)$ is
\textbf{$(e,i)$-isotropized} for some $i \in e \in\binom{[d]}{\geq
2}$ if
\[\Phi_e\cap \Big(\bigvee_{j\in[d]\setminus e}\Phi_{\{i,j\}}\Big) = \bigvee_{j\in[d]\setminus e}\Phi_{e\cup\{j\}}\]
up to $\mu$-negligible sets.

It is \textbf{fully isotropized} if it is $(e,i)$-isotropized for
every $(e,i)$.
\end{dfn}

Intuitively, both pleasantness and isotropizedness (say when $e =
[d]$) assert that the factors $\Phi_{\langle i\rangle}$ are `large
enough': in the first case, large enough to account for all of the
possible correlations between the coordinate projections under the
Furstenberg self-joining, and in the second to account for all of
the possible intersection between $\Phi_e$ and the combination
$\bigvee_{j\in[d]\setminus e}\Phi_{\{i,j\}}$ up to negligible sets.
This notion of pleasantness is very similar to Definition 4.2
in~\cite{Aus--nonconv}, where `pleasant systems' were first
introduced as those in which the larger factors $\S^{T_i}\vee
\Phi_{\langle i\rangle}$ were `characteristic' for the asymptotic
behaviour of the nonconventional averages~(\ref{eq:nonconv}) in
$L^2(\mu)$. Here our emphasis is rather different, since we are
concerned only with the integrals of these ergodic averages, rather
than the functions themselves. For these integrals it turns out that
we can discard the factors $\S^{T_i}$ from consideration. This
lightens some of the notation that follows, but otherwise makes very
little difference to the work we must go through.

Notice that the subset $e \subseteq [d]$ is allowed to vary in both
of the above definitions: this nuance is important, since the
pleasantness property relating a proper subfamily of actions $T_i$,
$i \in e$, is in general not a consequence of the pleasantness of
the whole family, and similarly for isotropizedness.

The main goal of this section is the following proposition.

\begin{prop}[Simultaneously pleasant and isotropized
extensions]\label{prop:pleasant-and-isotropized} Any commuting tuple
of actions $T_1,T_2,\ldots,T_d:\bbZ\curvearrowright (X,\S,\mu)$
admits an extension that is both fully pleasant and fully
isotropized.
\end{prop}

This will rely on a number of simpler steps, many closely following
the arguments of~\cite{Aus--nonconv}. We first show that any tuple
of actions admits an $(e,i)$-pleasant extension and, separately, an
$(e,i)$-isotropized extension.

The first of these results is proved exactly as was Proposition 4.6
in~\cite{Aus--nonconv}, and so we shall only sketch the proof here.
The idea behind the proof is to construct of a tower of extensions,
each accounting for the shortfall from pleasantness of its
predecessor, and then the pass to the inverse limit.

\begin{lem}[Existence of an $(e,i)$-pleasant
extension]\label{lem:partpleasant} Any commuting tuple of actions
$T_1,T_2,\ldots,T_d:\bbZ\curvearrowright (X,\S,\mu)$ admits an
$(e,i)$-pleasant extension
$(\tilde{X},\tilde{\S},\tilde{\mu},\tilde{T})$.
\end{lem}

\textbf{Proof}\quad We form
$(\tilde{X},\tilde{\S},\tilde{\mu},\tilde{T})$ as the inverse limit
of a tower of smaller extensions, each constructed from the
Furstenberg self-joining of its predecessor.  Let
$(X^{(1)},\S^{(1)},\mu^{(1)})$ be the Furstenberg self-joining
$(X^e,\S^{\otimes e},\mu^{\rm{F}}_e)$ and define on it the
transformations
\[\tilde{T}_i:=\prod_{j \in e}T_j\]
and
\[\tilde{T}_k:=(T_k)^{\times e}\quad\quad\hbox{for}\ k \neq i,\]
and interpret it as an extension of $(X,\S,\mu,T)$ with the
coordinate projection $\pi_i$ as factor map. We now see that if $f_j
\in L^\infty(\mu)$ for each $j \in e$ then
\[\int_{X^e}\prod_{j\in e}f_j\circ\pi_j\,\d\mu^{\rm{F}}_e = \int_{X^e}\sfE_\mu\big(f_i\circ\pi_i\,\big|\,(\pi_j)_{j\in e\setminus\{i\}}\big)\cdot\prod_{j\in
e\setminus\{i\}}f_j\circ\pi_j\,\d\mu^{\rm{F}}_e,\] and from the
above definition that the factor of $X^e = X^{(1)}$ generated by
$(\pi_j)_{j\in e\setminus\{i\}}$ is contained in $\bigvee_{j\in
e\setminus \{i\}}\Phi^{(1)}_{\{i,j\}}$. If we now iterate this
construction to form $(X^{(2)},\S^{(2)},\mu^{(2)},T^{(2)})$ from
$(X^{(1)},\S^{(1)},\mu^{(1)},T^{(1)})$, and so on, then the
approximation argument given for Proposition 4.6
of~\cite{Aus--nonconv} shows that the inverse limit is
$(e,i)$-pleasant. \qed

\textbf{Remark}\quad Since the appearance of~\cite{Aus--nonconv},
Bernard Host has given in~\cite{Hos??} a method for constructing a
pleasant extension of a system without recourse to an inverse limit.
However, we will make further use of inverse limits momentarily to
construct an extension that is fully pleasant, rather than just
$(e,i)$-pleasant for some fixed $(e,i)$, and at present we do not
know of any quicker construction guaranteeing this stronger
condition. \fin

A similar argument gives the existence of $(e,i)$-isotropized
extensions.

\begin{lem}[Existence of $(e,i)$-isotropized
extension]\label{lem:isotropized} Any commuting tuple of actions
$T_1,T_2,\ldots,T_d:\bbZ\curvearrowright (X,\S,\mu)$ admits an
$(e,i)$-isotropized extension
$(\tilde{X},\tilde{\S},\tilde{\mu},\tilde{T})$.
\end{lem}

\textbf{Proof}\quad Once again we build this as an inverse limit.
First form the relatively independent self-product
$(X^{(1)},\mu^{(1)}):=(X^2,\S\otimes_{\Phi_e}\S,\mu\otimes_{\Phi_e}\mu)$
with coordinate projections $\pi_1$, $\pi_2$ back onto $(X,\S,\mu)$,
and interpret it as an extension of $(X,\S,\mu)$ through the first
of these.  Choose arbitrarily some $i \in e$, and now define the
extended actions $T^{(1)}_j$ on $X^{(1)}$ by setting \[T^{(1)}_j
:=\left\{\begin{array}{ll}T_j\times T_j&\quad\hbox{if}\ j\not\in e,\\
T_j\times T_i&\quad\hbox{if}\ j\in e;\end{array}\right.\] these all
preserve $\mu^{(1)}$, even in the latter case, because our product
is relatively independent over the factor left invariant by each
$T_j^{-1}T_i$ for $j \in e$.

We now extend $(X^{(1)},\S^{(1)},\mu^{(1)},T^{(1)})$ to
$(X^{(2)},\S^{(2)},\mu^{(2)},T^{(2)})$ by repeating the same
construction, and so on, to form an inverse series with inverse
limit $(\tilde{X},\tilde{\S},\tilde{\mu},\tilde{T})$.

We will show that this has the desired property. Any $f \in
{L^\infty(\mu\uhr_{\Phi_e\cap (\bigvee_{j\in[d]\setminus
e}\Phi_{\{i,j\}})})}$ may, in particular, be approximated in
$L^1(\mu)$ by finite sums of the form $\sum_p\prod_{j\in[d]\setminus
e}\phi_{p,j}$ with $\phi_{j,p} \in
L^\infty(\mu\uhr_{\Phi_{\{i,j\}}})$.  However, since $\mu^{(1)}$ is
joined relatively independently conditioned on $\Phi_e$ and $f$ is
also $\Phi_e$-measurable, it follows that $f\circ\pi_1=f\circ\pi_2$
$\mu^{(1)}$-almost surely, and so in the extended system
$(X^{(1)},\S^{(1)},\mu^{(1)})$ we can alternatively approximate
$f\circ\pi_1$ by the functions $\sum_p\prod_{j\in[d]\setminus
e}\phi_{p,j}\circ\pi_2$; and now every $\phi_{p,j}\circ\pi_2$ is
both manifestly $\Phi^{(1)}_{\{j,i\}}$-measurable, since both $T_j$
and $T_i$ are simply lifted to $T_j^{\times 2}$ and $T_i^{\times
2}$, and manifestly $\Phi^{(1)}_e$-measurable, since all the
transformations $T^{(1)}_j$ defined above for $j\in e$ agree on the
second coordinate factor $\pi_2^{-1}(\S)$. Therefore $f\circ\pi_1$
may be approximated arbitrarily well in $L^1(\mu)$ by functions that
are measurable with respect to $\bigvee_{j\in[d]\setminus
e}\Phi^{(1)}_{e\cup\{j\}}$. Now another simple approximation
argument and the martingale convergence theorem show that the
inverse limit system $(\tilde{X},\tilde{\S},\tilde{\mu},\tilde{T})$
is actually $(e,i)$-isotropized, as required. \qed

We will finish the proof of
Proposition~\ref{prop:pleasant-and-isotropized} using the following
properties of stability under forming further inverse limits.

\begin{lem}[Pleasantness of inverse
limits]\label{lem:pleasant-inv-lim} If
\[\ldots\to(X^{(m+1)},\S^{(m+1)},\mu^{(m+1)},T^{(m+1)})\stackrel{\psi_m}{\longrightarrow}(X^{(m)},\S^{(m)},\mu^{(m)},T^{(m)})\to\ldots\]
is an inverse system with inverse limit
$(\tilde{X},\tilde{\S},\tilde{\mu},\tilde{T})$ and $i \in e
\subseteq [d]$, then
\begin{itemize}
\item if $(X^{(m)},\S^{(m)},\mu^{(m)},T^{(m)})$ is $(e,i)$-pleasant for
infinitely many $m$, then
$(\tilde{X},\tilde{\S},\tilde{\mu},\tilde{T})$ is also
$(e,i)$-pleasant;
\item if $(X^{(m)},\S^{(m)},\mu^{(m)},T^{(m)})$ is $(e,i)$-isotropized for
infinitely many $m$, then
$(\tilde{X},\tilde{\S},\tilde{\mu},\tilde{T})$ is also
$(e,i)$-isotropized.
\end{itemize}
\end{lem}

\textbf{Proof}\quad We give the proof for the retention of
$(e,i)$-pleasantness, the case of $(e,i)$-isotropizedness being
exactly analogous.

Since any $1$-bounded member of $L^\infty(\tilde{\mu})$ may be
approximated arbitrarily well in $L^1(\tilde{\mu})$ by $1$-bounded
members of $L^\infty(\mu^{(m)})$, by a simple approximation argument
it will suffice to prove that given $m\geq 1$ and $f_j \in
L^\infty(\mu^{(m)})$ for each $j\in e$ we have
\[\int_X \prod_{j\in e}f_j\circ\pi_j\,\d\tilde{\mu}^{\rm{F}}_e = \int_X \sfE_{\tilde{\mu}}\Big(f_i\,\Big|\,\bigvee_{j\in e\setminus\{i\}}\tilde{\Phi}_{\{i,j\}}\Big)\circ\pi_i\cdot \prod_{j\neq i}f_j\circ\pi_j\,\d\tilde{\mu}^{\rm{F}}_e.\]
However, by definition and Lemma~\ref{lem:Fberg-inv-lim} we know
that after choosing any $m_1 \geq m$ for which
$(X^{(m_1)},\S^{(m_1)},\mu^{(m_1)},T^{(m_1)})$ is $(e,i)$-pleasant
the above is obtained with $\bigvee_{j\in
e\setminus\{i\}}\Phi^{(m_1)}_{\{i,j\}}$ in place of $\bigvee_{j\in
e\setminus\{i\}}\tilde{\Phi}_{\{i,j\}}$, and now letting $m_1 \to
\infty$ and appealing to the bounded martingale convergence theorem
gives the result. \qed

It now remains only to collect our different properties together
using more inverse limits, whose organization is now rather
arbitrary.

\textbf{Proof of
Proposition~\ref{prop:pleasant-and-isotropized}}\quad Pick a
sequence of pairs $((e_m,i_m))_{m\geq 1}$ from the finite set
$\{(e,i):\ |e|\geq 2,\,i \in e \subseteq [d]\}$ in which each
possible $(e,i)$ appears infinitely often. Now one last time form a
tower of extensions
\[\ldots\to(X^{(m+1)},\S^{(m+1)},\mu^{(m+1)},T^{(m+1)})\to (X^{(m)},\S^{(m)},\mu^{(m)},T^{(m)})\to\ldots\]
above $(X,\S,\mu,T)$ in which $(X^{(m)},\S^{(m)},\mu^{(m)},T^{(m)})$
is $(e_{(m+1)/2},i_{(m+1)/2})$-pleasant when $m$ is odd and
$(e_{m/2},i_{m/2})$-isotropized when $m$ is even.  The two parts of
Lemma~\ref{lem:pleasant-inv-lim} now show that the resulting inverse
limit extension has all the desired properties. \qed

\section{Furstenberg self-joinings of pleasant and isotropized systems}

Having established that all systems have fully pleasant and
isotropized extensions, it remains to explain the usefulness of such
extensions for the proof of Theorem~\ref{thm:multirec}.  This
derives from the implications of these conditions for the structure
of the Furstenberg self-joining.

\begin{lem}\label{lem:Fberg-struct}
If the tuple $T_1,T_2,\ldots,T_d:\bbZ\curvearrowright(X,\S,\mu)$ is
fully pleasant and fully isotropized,
$\cal{I}\subseteq\binom{[d]}{\geq 2}$ is an up-set and $e$ is a
member of $\binom{[d]}{\geq 2}\setminus\cal{I}$ of maximal size then
the oblique copy $\Phi^{\rm{F}}_e$ and the oblique factor
$\Phi^{\rm{F}}_{\cal{I}}$ are relatively independent over
$\Phi^{\rm{F}}_{\cal{I}\cap\langle e\rangle}$ under
$\mu^{\rm{F}}_{[d]}$.
\end{lem}

\textbf{Proof}\quad Suppose that $F_1 \in
L^\infty(\mu^{\rm{F}}_{[d]}\uhr_{\Phi^{\rm{F}}_e})$ and $F_2 \in
L^\infty(\mu^{\rm{F}}_{[d]}\uhr_{\Phi^{\rm{F}}_{\cal{I}}})$. It will
suffice to show that
\[\int_{X^d}F_1F_2\,\d\mu^{\rm{F}}_{[d]} = \int_{X^d}\sfE_\mu(F_1\,|\,\Phi^{\rm{F}}_{\cal{I}\cap\langle
e\rangle})\cdot F_2\,\d\mu^{\rm{F}}_{[d]}.\]

Pick $i\in e$. By Lemma~\ref{lem:diag} there is some $f_1 \in
L^\infty(\mu\uhr_{\Phi_e})$ such that $F_1 = f_1\circ\pi_i$
$\mu^{\rm{F}}_{[d]}$-almost surely.

Let $\{a_1,a_2,\ldots,a_k\}$ be the antichain of minimal elements in
$\cal{I}$; this clearly generates $\cal{I}$ as an up-set. Since $e
\not\in\cal{I}$ we must have $a_j \setminus e\not=\emptyset$ for
each $j \leq k$. Pick $i_j \in a_j\setminus e$ arbitrarily for each
$j\leq k$, so that, again by Lemma~\ref{lem:diag},
$\Phi^{\rm{F}}_{a_j} = \pi_{i_j}^{-1}(\Phi_{a_j})$ (up to
$\mu^{\rm{F}}_{[d]}$-negligible sets).

Now, since $\Phi^{\rm{F}}_{\cal{I}} = \bigvee_{j\leq
k}\Phi^{\rm{F}}_{a_j}$, $F_2$ may be approximated arbitrarily well
in $L^1(\mu^{\rm{F}}_{[d]})$ by sums of products of the form
$\sum_p\prod_{j\leq k}\phi_{j,p}\circ\pi_{i_j}$ with $\phi_{j,p}\in
L^\infty(\mu\uhr_{\Phi_{a_j}})$, and so by continuity and linearity
it suffices to assume that $F_2$ is an individual such product term.
This represents $F_2$ as a function of coordinates in $X^d$ indexed
only by members of $[d]\setminus e$, and now we appeal to
Lemma~\ref{lem:Fberg-project} and the pleasantness of
$\mu^{\rm{F}}_{([d]\setminus e)\cup\{i\}}$ to deduce that
\begin{multline*}
\int_{X^d}F_1\cdot \prod_{j\leq
k}\phi_{j,p}\circ\pi_{i_j}\,\d\mu^{\rm{F}}_{[d]}\\ =
\int_{X^{([d]\setminus
e)\cup\{i\}}}\sfE_\mu\Big(f_1\,\Big|\,\bigvee_{j\in[d]\setminus
e}\Phi_{\{i,j\}}\Big)\cdot \prod_{j\leq
k}\phi_{j,p}\circ\pi_{i_j}\,\d\mu^{\rm{F}}_{([d]\setminus
e)\cup\{i\}}.
\end{multline*}

However, now the property that $(X,\S,\mu,T)$ is $(e,i)$-isotropized
and the fact that $f_1$ is already $\Phi_e$-measurable imply that
\[\sfE_\mu\Big(f_1\,\Big|\,\bigvee_{j\in[d]\setminus e}\Phi_{\{i,j\}}\Big) = \sfE_\mu\Big(f_1\,\Big|\,\bigvee_{j\in [d]\setminus e} \Phi_{e \cup \{j\}}\Big),\]
and since each $e \cup \{j\} \in \cal{I}$ (by the maximality of $e$
in $\cal{P}[d]\setminus\cal{I}$), under $\pi_i$ this conditional
expectation must be identified with
$\sfE_\mu(F_1\,|\,\Phi^{\rm{F}}_{\cal{I}\cap\langle e\rangle})$, as
required. \qed

\begin{cor}\label{cor:Fberg-struct}
If the tuple $T_1,T_2,\ldots,T_d:\bbZ\curvearrowright(X,\S,\mu)$ is
fully pleasant and fully isotropized and $\cal{I},\cal{I}'\subseteq
\binom{[d]}{\geq 2}$ are two up-sets then $\Phi^{\rm{F}}_{\cal{I}}$
and $\Phi^{\rm{F}}_{\cal{I}'}$ are relatively independent over
$\Phi^{\rm{F}}_{\cal{I}\cap\cal{I}'}$ under $\mu^{\rm{F}}_{[d]}$.
\end{cor}

\textbf{Proof}\quad This is proved for fixed $\cal{I}$ by induction
on $\cal{I}'$. If $\cal{I}' \subseteq \cal{I}$ then the result is
clear, so now let $e$ be a minimal member of
$\cal{I}'\setminus\cal{I}$ of maximal size, and let $\cal{I}'' :=
\cal{I}'\setminus\{e\}$.  It will suffice to prove that if $F \in
L^\infty(\mu^{\rm{F}}_{[d]}\uhr_{\Phi^{\rm{F}}_{\cal{I}'}})$ then
\[\sfE_{\mu^{\rm{F}}_{[d]}}(F\,|\,\Phi^{\rm{F}}_{\cal{I}}) = \sfE_{\mu^{\rm{F}}_{[d]}}(F\,|\,\Phi^{\rm{F}}_{\cal{I}\cap\cal{I}'}),\]
and furthermore, by approximation, to do so only for $F$ that are of
the form $F_1\cdot F_2$ with $F_1 \in
L^\infty(\mu^{\rm{F}}_{[d]}\uhr_{\Phi^{\rm{F}}_{\langle e
\rangle}})$ and $F_2 \in
L^\infty(\mu^{\rm{F}}_{[d]}\uhr_{\Phi^{\rm{F}}_{\cal{I}''}})$.
However, for these we can write \begin{multline*}
\sfE_{\mu^{\rm{F}}_{[d]}}(F\,|\,\Phi^{\rm{F}}_{\cal{I}}) =
\sfE_{\mu^{\rm{F}}_{[d]}}\big(\sfE_{\mu^{\rm{F}}_{[d]}}(F\,|\,\Phi^{\rm{F}}_{\cal{I}\cup\cal{I}''})\,\big|\,\Phi^{\rm{F}}_{\cal{I}}\big)\\
=
\sfE_{\mu^{\rm{F}}_{[d]}}\big(\sfE_{\mu^{\rm{F}}_{[d]}}(F_1\,|\,\Phi^{\rm{F}}_{\cal{I}\cup\cal{I}''})\cdot
F_2\,\big|\,\Phi^{\rm{F}}_{\cal{I}}\big),
\end{multline*}
and by Lemma~\ref{lem:Fberg-struct}
\[\sfE_{\mu^{\rm{F}}_{[d]}}(F_1\,|\,\Phi^{\rm{F}}_{\cal{I}\cup\cal{I}''})
=
\sfE_{\mu^{\rm{F}}_{[d]}}(F_1\,|\,\Phi^{\rm{F}}_{(\cal{I}\cup\cal{I}'')\cap\langle
e\rangle}).\] On the other hand $(\cal{I}\cup\cal{I}'')\cap\langle
e\rangle \subseteq \cal{I}''$ (because $\cal{I}''$ contains every
subset of $[d]$ that strictly includes $e$, since $\I'$ is an
up-set), and so Lemma~\ref{lem:Fberg-struct} promises similarly that
\[\sfE_{\mu^{\rm{F}}_{[d]}}(F_1\,|\,\Phi^{\rm{F}}_{(\cal{I}\cup\cal{I}'')\cap\langle
e\rangle}) =
\sfE_{\mu^{\rm{F}}_{[d]}}(F_1\,|\,\Phi^{\rm{F}}_{\cal{I}''}).\]
Therefore the above expression for
$\sfE_{\mu^{\rm{F}}_{[d]}}(F\,|\,\Phi^{\rm{F}}_{\cal{I}})$
simplifies to
\begin{multline*}
\sfE_{\mu^{\rm{F}}_{[d]}}\big(\sfE_{\mu^{\rm{F}}_{[d]}}(F_1\,|\,\Phi^{\rm{F}}_{\cal{I}''})\cdot
F_2\,\big|\,\Phi^{\rm{F}}_{\cal{I}}\big) =
\sfE_{\mu^{\rm{F}}_{[d]}}\big(\sfE_{\mu^{\rm{F}}_{[d]}}(F_1\cdot
F_2\,|\,\Phi^{\rm{F}}_{\cal{I}''})\,\big|\,\Phi^{\rm{F}}_{\cal{I}}\big)\\
=
\sfE_{\mu^{\rm{F}}_{[d]}}\big(\sfE_{\mu^{\rm{F}}_{[d]}}(F\,|\,\Phi^{\rm{F}}_{\cal{I}''})\,\big|\,\Phi^{\rm{F}}_{\cal{I}}\big)
= \sfE_{\mu^{\rm{F}}_{[d]}}(F\,|\,\Phi^{\rm{F}}_{\cal{I}\cap
\cal{I}''}) =
\sfE_{\mu^{\rm{F}}_{[d]}}(F\,|\,\Phi^{\rm{F}}_{\cal{I}\cap\cal{I}'}),
\end{multline*}
by the inductive hypothesis applied to $\cal{I}''$ and $\cal{I}$, as
required. \qed

\section{Completion of the proof}

We have now set the stage for our analog of Tao's infinitary
hypergraph removal machinery.  Observe first that the conclusion of
Theorem~\ref{thm:multirec} clearly holds for the commuting tuple
$T_1,T_2,\ldots,T_d:\bbZ\curvearrowright (X,\S,\mu)$ if it holds for
any extension of that tuple.  Therefore by
Proposition~\ref{prop:pleasant-and-isotropized} we may assume our
commuting tuple is fully pleasant and fully isotropized, and so need
only prove for such $\mu$ that if $\mu (A)
> 0$ then $\mu^{\rm{F}}_{[d]}(A^d) > 0$. For these
particular $\mu$ Corollary~\ref{cor:Fberg-struct} gives us a very
precise picture of the joint law under $\mu^{\rm{F}}_{[d]}$ of the
poset of oblique factors $\Phi^{\rm{F}}_{\I}$, and hence actually of
the inverse image factors $\pi_i^{-1}(\Phi_{\I})$ for $\I\subseteq
\langle i\rangle$.

Note that we have not tamed all of the potentially wild structure of
the joint distribution of the factors $\Phi_\cal{I}$ under $\mu$,
but only that of the associated oblique factors under the
Furstenberg self-joining $\mu^{\rm{F}}_{[d]}$. It seems quite likely
that in some cases the factors $\Phi_\I$ of the original system can
still exhibit a very complicated joint distribution, even after
passing to a fully pleasant and isotropized extension. However, the
understanding of the oblique copies is already enough to complete
the proof of multiple recurrence using a relative of Tao's
`infinitary removal lemma' in~\cite{Tao07}. One of his chief
innovations was an infinitary analog of the property of hypergraph
removability for a collection of factors of a probability space
(Theorem 4.2 of~\cite{Tao07}). Here we shall actually make do with a
more modest conclusion than his `removability', but our argument
will follow essentially the same steps. We shall derive
Theorem~\ref{thm:multirec} as the top case of the following
inductive claim, tailored to our present needs.

\begin{prop}\label{prop:infremoval}
Suppose that $\I_{i,j}$ for $i=1,2,\ldots,d$ and $j =
1,2,\ldots,k_i$ are collections of up-sets in $\binom{[d]}{\geq 2}$
such that $[d] \in \I_{i,j} \subseteq \langle i\rangle$ for each
$i,j$, and suppose further that the sets $A_{i,j}\in
\Phi_{\I_{i,j}}$ are such that
\[\mu^{\rm{F}}_{[d]}\Big(\prod_{i=1}^d\Big(\bigcap_{j=1}^{k_i}A_{i,j}\Big)\Big) = 0.\]
Then we must also have
\[\mu\Big(\bigcap_{i=1}^d\bigcap_{j=1}^{k_i}A_{i,j}\Big) = 0.\]
\end{prop}

The following terminology will be convenient during the proof.

\begin{dfn}
A family $(\I_{i,j})_{i,j}$ has the property \textbf{P} if it
satisfies the conclusion of the preceding proposition.
\end{dfn}

The conclusion of multiple recurrence follows from
Proposition~\ref{prop:infremoval} at once:

\textbf{Proof of Theorem~\ref{thm:multirec} from
Proposition~\ref{prop:infremoval}}\quad Suppose that $A \in \S$ is
such that $\mu^{\rm{F}}_{[d]}(A^d) = 0$. Then by the pleasantness of
the whole system we have
\[\mu^{\rm{F}}_{[d]}(A^d) = \int_{X^d}\prod_{i=1}^d \sfE_{\mu}(1_A\,|\,\Phi_{\langle i\rangle})\circ\pi_i\,\d\mu^{\rm{F}}_{[d]} = 0.\]
Now the level set $B_i := \{\sfE_{\mu}(1_A\,|\,\Phi_{\langle
i\rangle})
> 0\}$ (of course, this is unique only up to $\mu$-negligible sets) lies in $\Phi_{\langle i\rangle}$, and the above equality certainly
implies that also $\mu^{\rm{F}}_{[d]}(B_1\times B_2\times
\cdots\times B_d) = 0$.  Now, on the one hand, setting $k_i=1$,
$\I_{i,1}:= \langle i\rangle$ and $A_{i,1} := B_i$ for each $i \leq
d$, Proposition~\ref{prop:infremoval} tells us that $\mu(B_1\cap
B_2\cap\cdots\cap B_d) = 0$, while on the other we must have
$\mu(A\setminus B_i) = 0$ for each $i$, and so overall $\mu(A) \leq
\mu(B_1\cap B_2\cap\cdots\cap B_d) + \sum_{i=1}^d\mu(A\setminus B_i)
= 0$, as required. \qed

It remains to prove Proposition~\ref{prop:infremoval}.  This will be
done by induction on a suitable ordering of the possible collections
of up-sets $(\I_{i,j})_{i,j}$, appealing to a handful of different
possible cases at different steps of the induction.  At the
outermost level, this induction will be organized according to the
depth of our up-sets (defined in Section~\ref{sec:prelim}).

Let us first illustrate how the above reduction to
Proposition~\ref{prop:infremoval} and then the inductive proof of
that proposition combine to give a proof of
Theorem~\ref{thm:multirec} in the simple case $d=3$.

\textbf{Example}\quad Suppose that $T_1,T_2,T_3:\bbZ\curvearrowright
(X,\S,\mu)$ is a fully pleasant and fully isotropized triple of
actions and that $A \in \S$ has $\mu^{\rm{F}}_{[3]}(A^3) = 0$.  We
will show that $\mu(A) = 0$.  As in the above argument, we know that
\[\mu^{\rm{F}}_{[3]}(A^3) = \int_{X^3}(\sfE_\mu(1_A\,|\,\Phi_{\langle 1\rangle})\circ\pi_1)\cdot(\sfE_\mu(1_A\,|\,\Phi_{\langle 2\rangle})\circ\pi_2)\cdot(\sfE_\mu(1_A\,|\,\Phi_{\langle 3\rangle})\circ\pi_3)\,\d\mu^{\rm{F}}_{[3]},\]
and so we must actually have $\mu^{\rm{F}}_{[d]}\big(B_1\times
B_2\times B_3) = 0$ where $B_i:= \{\sfE_\mu(1_A\,|\,\Phi_{\langle
i\rangle}) > 0\}$. Clearly $A$ is contained in $B_1\cap B_2\cap B_3$
up to a $\mu$-negligible set, so it will suffice to show that this
intersection is $\mu$-negligible.

Now, each of $\Phi_{\langle 1\rangle} =
\Phi_{\{1,2\}}\vee\Phi_{\{1,3\}}$, $\Phi_{\langle 2\rangle}$ and
$\Phi_{\langle 3\rangle}$ can be generated using intersections of
members from countable generating sets in each $\Phi_{\{i,j\}}$. Let
\[\B_{\{i,j\},1}\subseteq \B_{\{i,j\},2}\subseteq \ldots\]
be an increasing sequence of finite subalgebras of sets that
generates $\Phi_{\{i,j\}}$ up to $\mu$-negligible sets, and let
\[\Xi^{(n)}_i := \S^{T_1 = T_2 = T_3}\vee\B_{\{i,j\},n}\vee\B_{\{i,k\},n}\]
when $\{i,j,k\} = \{1,2,3\}$. By the martingale convergence theorem
we have $\sfE_\mu(1_{B_i}\,|\,\Xi^{(n)}_i) \to 1_{B_i}$ in
$L^2(\mu)$ as $n\to\infty$.  Now pick $\delta \in (0,1/3)$ and let
$C_i^{(n)} := \{\sfE_\mu(1_{B_i}\,|\,\Xi^{(n)}_i) > 1 - \delta\}$,
so for large $n$ this set should be a $\mu$-approximation to $B_i$,
and observe in addition that $\sfE_\mu(1_{C_1^{(n)}\setminus
B_1}\,|\,\Xi^{(n)}_1)\leq \delta$ almost surely.

It easy to check from Corollary~\ref{cor:Fberg-struct} that
$\Phi_{\langle i\rangle}^{\rm{F}}$ must be relatively independent
from $\pi_i^{-1}(\Xi^{(n)}_j\vee\Xi^{(n)}_k)$ over
$\pi_i^{-1}(\Xi^{(n)}_i)$ under $\mu^{\rm{F}}_{[3]}$ when $\{i,j,k\}
= \{1,2,3\}$, and from this we compute that
\begin{multline*}
\mu^{\rm{F}}_{[3]}((C_1^{(n)}\times C_2^{(n)}\times
C_3^{(n)})\setminus \pi_1^{-1}(B_1)) =
\int_{X^3}(\sfE_\mu(1_{C_1^{(n)}\setminus
B_1}\,|\,\Xi^{(n)}_1)\circ\pi_1)\cdot 1_{C^{(n)}_2}\cdot
1_{C^{(n)}_3}\,\d\mu^{\rm{F}}_{[3]}\\
\leq \delta\int_{X^3}1_{C^{(n)}_1}\cdot 1_{C^{(n)}_2}\cdot
1_{C^{(n)}_3}\,\d\mu^{\rm{F}}_{[3]} =
\delta\mu^{\rm{F}}_{[3]}(C_1^{(n)}\times C^{(n)}_2\times C^{(n)}_3).
\end{multline*}
Therefore
\[\mu^{\rm{F}}_{[3]}(C_1^{(n)}\times C^{(n)}_2\times C^{(n)}_3) \leq \mu^{\rm{F}}_{[d]}(B_1 \times B_2\times B_3) + 3\delta\mu^{\rm{F}}_{[d]}(C_1^{(n)}\times C^{(n)}_2\times C^{(n)}_3),\]
and so since $\delta < 1/3$ we must have
$\mu^{\rm{F}}_{[3]}(C_1^{(n)}\times C^{(n)}_2\times C^{(n)}_3) = 0$
for all $n$.

The importance of this is that for large $n$ we have now
approximated the sets $B_i$ by sets $C_i^{(n)}$ that lie in the
simpler $\s$-algebras $\Xi^{(n)}_i$ but nevertheless still enjoy the
property that the measure $\mu^{\rm{F}}_{[3]}(C_1^{(n)}\times
C^{(n)}_2\times C^{(n)}_3)$ is \emph{strictly} zero. Since each
$\B_{\{i,j\},n}$ is finite, for any given $n$ we may write each
$C_i^{(n)}$ as a finite union of subsets of the form $C_{i,p}^{(n)}
= D_{i,p}\cap C_{i,j,p}\cap C_{i,k,p}$ with $D_{i,p} \in \S^{T_1 =
T_2 = T_3}$ and $C_{i,j,p} \in \B_{\{i,j\},n}$ for every $p$, and
these must now also enjoy the property that
\[\mu^{\rm{F}}_{[3]}\big((D_{1,p_1}\cap C_{1,2,p_1}\cap C_{1,3,p_1})\times (D_{1,p_2}\cap C_{2,1,p_2}\cap C_{2,3,p_2})\times (D_{3,p_3}\cap C_{3,1,p_3}\cap C_{3,2,p_3})\big) = 0\]
for all possible indices $p_1$, $p_2$, $p_3$.

Next the fact that
$\mu^{\rm{F}}_{[3]}(\pi_i^{-1}(C)\triangle\pi_j^{-1}(C)) = 0$
whenever $C \in \Phi_{\{i,j\}}$ (Lemma~\ref{lem:diag}) comes into
play, allowing us for example to move the set $C_{2,1,p_2}$ under
the first coordinate rather than the second in the above equation,
and similarly. In this way we can re-arrange the above equation into
the form
\begin{multline*}
\mu^{\rm{F}}_{[3]}\big(((D_{1,p_1}\cap D_{1,p_2}\cap D_{1,p_3})\cap
(C_{1,2,p_1}\cap C_{2,1,p_2})\cap (C_{1,3,p_1}\cap C_{3,1,p_3}))\\
\times (C_{2,3,p_2}\cap C_{3,2,p_3})\times X\big) = 0.
\end{multline*}

This equation involves the sets $D := D_{1,p_1}\cap D_{1,p_2}\cap
D_{1,p_3} \in \S^{T_1 = T_2 = T_3}$, $C_{1,2} := C_{1,2,p_1}\cap
C_{2,1,p_2} \in \Phi_{\{1,2\}}$, $C_{1,3} := C_{1,3,p_1}\cap
C_{3,1,p_3} \in \Phi_{\{1,3\}}$ and $C_{2,3} := C_{2,3,p_2}\cap
C_{3,2,p_3} \in \Phi_{\{2,3\}}$.  Now,
Corollary~\ref{cor:Fberg-struct} tells us that the three oblique
copies $\Phi_{\{i,j\}}^{\rm{F}}$ are relatively independent over
$\Phi_{\{1,2,3\}}^{\rm{F}}$ under $\mu^{\rm{F}}_{[3]}$, and so we
deduce from the above equation that
\begin{eqnarray*}
0&=&\int_{X^3}(1_D\circ\pi_1)\cdot(\sfE_\mu(1_{C_{1,2}}\,|\,\S^{T_1
=
T_2 = T_3})\circ\pi_1)\\
&&\quad\quad\quad\quad\quad\quad\cdot(\sfE_\mu(1_{C_{1,3}}\,|\,\S^{T_1
= T_2 = T_3})\circ\pi_1)\cdot (\sfE_\mu(1_{C_{2,3}}\,|\,\S^{T_1 =
T_2 =
T_3})\circ\pi_2)\,\d\mu^{\rm{F}}_{[3]}\\
&=& \int_X 1_D\cdot\sfE_\mu(1_{C_{1,2}}\,|\,\S^{T_1 = T_2 =
T_3})\cdot\sfE_\mu(1_{C_{1,3}}\,|\,\S^{T_1 = T_2 = T_3})\cdot
\sfE_\mu(1_{C_{2,3}}\,|\,\S^{T_1 = T_2 = T_3})\,\d\mu
\end{eqnarray*}
where the first and second line here are equal by
Lemma~\ref{lem:diag} since all the functions involved are $\S^{T_1 =
T_2 = T_3}$-measurable.

However, this now implies that
\[\mu\Big(D\cap \bigcap_{\{i,j\}\in\binom{[3]}{2}}\{\sfE_\mu(1_{C_{i,j}}\,|\,\S^{T_1 = T_2 =
T_3}) > 0\}\Big) = 0\] and hence that we must also have
\[\mu\big(D_{1,p_1}\cap D_{1,p_2}\cap
D_{1,p_3}\cap C_{1,2,p_1}\cap C_{2,1,p_2}\cap C_{1,3,p_1}\cap
C_{3,1,p_3}\cap C_{2,3,p_2}\cap C_{3,2,p_3}\big) = 0.\] Taking the
union of these equations over triples of indices $p_1$, $p_2$, $p_3$
gives $\mu(C_1^{(n)}\cap C_2^{(n)}\cap C_3^{(n)}) = 0$ for any $n$,
and so since the sets $C_i^{(n)}$ approximate $B_i$ as $n\to\infty$
it follows that $\mu(B_1\cap B_2\cap B_3) = 0$, as required. \fin

We now turn to full induction that generalizes the above argument,
broken into a number of steps.

\begin{lem}[Lifting using relative independence]\label{lem:lift-rel-ind}
Suppose that all up-sets in the collection $(\I_{i,j})_{i,j}$ have
depth at least $k$, that all those with depth exactly $k$ are
principal, and that there are $\ell \geq 1$ of these.  Then if
property P holds for all similar collections having $\ell - 1$
up-sets of depth $k$, then it holds also for this collection.
\end{lem}

\textbf{Proof}\quad Let $\I_{i_1,j_1} = \langle e_1\rangle$,
$\I_{i_2,j_2} = \langle e_2\rangle$, \ldots, $\I_{i_\ell,j_\ell} =
\langle e_\ell\rangle$ be an enumeration of all the (principal)
up-sets of depth $k$ in our collection.  We will treat two separate
cases.

First suppose that two of the generating sets agree; by re-ordering
if necessary we may assume that $e_1 = e_2$.  Clearly we can assume
that there are no duplicates among the coordinate-collections
$(\I_{i,j})_{j=1}^{k_i}$ for each $i$ separately, so we must have
$i_1 \neq i_2$. However, if we now suppose that $A_{i,j}\in
\I_{i,j}$ for each $i$, $j$ are such that
\[\mu^{\rm{F}}_{[d]}\Big(\prod_{i=1}^d\Big(\bigcap_{j=1}^{k_i}A_{i,j}\Big)\Big) = 0,\]
then the same equality holds if we simply replace $A_{i_1,j_1} \in
\langle e_1\rangle$ with $A'_{i_1,j_1}:= A_{i_1,j_1}\cap
A_{i_2,j_2}$ and $A_{i_2,j_2}$ with $A'_{i_2,j_2} := X$. Now this
last set can simply be ignored to leave an instance of a
$\mu^{\rm{F}}_{[d]}$-negligible product for the same collection of
up-sets omitting $\I_{i_2,j_2}$, and so property P of this reduced
collection completes the proof.

On the other hand, if all the $e_i$ are distinct, we shall simplify
the last of the principal up-sets $\I_{i_\ell,j_\ell}$ by exploiting
the relative independence among the associated oblique copies of our
factors. Assume for notational simplicity that $(i_\ell,j_\ell) =
(1,1)$; clearly this will not affect the proof.  We will reduce to
an instance of property P associated to the collection $(\I'_{i,j})$
defined by
\[\I'_{i,j} := \left\{\begin{array}{ll}\langle e_\ell\rangle\setminus\{e_\ell\}&\quad\hbox{if}\ (i,j) = (1,1)\\ \I_{i,j}&\quad\hbox{else,}\end{array}\right.\]
which has one fewer up-set of depth $k$ and so falls under the
inductive assumption.

Indeeed, we know from Corollary~\ref{cor:Fberg-struct} that under
$\mu^{\rm{F}}_{[d]}$ the set $\pi_1^{-1}(A_{1,1})$ is relatively
independent from all the sets $\pi_i^{-1}(A_{i,j})$, $(i,j) \neq
(1,1)$, over the factor $\pi_1^{-1}(\Phi_{\langle
e_\ell\rangle\setminus\{e_\ell\}})$, which is dense inside the
relevant oblique copy $\Phi^{\rm{F}}_{\langle
e_\ell\rangle\setminus\{e_\ell\}}$.  Therefore
\begin{multline*}
0 =
\mu^{\rm{F}}_{[d]}\Big(\prod_{i=1}^d\Big(\bigcap_{j=1}^{k_i}A_{i,j}\Big)\Big)\\
= \int_{X^d}\sfE_\mu(1_{A_{1,1}}\,|\,\Phi_{\langle
e_\ell\rangle\setminus\{e_\ell\}})\circ\pi_1\cdot\prod_{j=2}^{k_1}1_{\pi_1^{-1}(A_{1,j})}\cdot\prod_{i=2}^d\prod_{j=1}^{k_i}1_{\pi_i^{-1}(A_{i,j})}\,\d\mu^{\rm{F}}_{[d]}.
\end{multline*}
Setting $A'_{1,1}:= \{\sfE_\mu(1_{A_{1,1}}\,|\,\Phi_{\langle
e_\ell\rangle\setminus\{e_\ell\}}) > 0\} \in \Phi_{\langle
e_\ell\rangle\setminus\{e_\ell\}}$ and $A'_{i,j} := A_{i,j}$ for
$(i,j) \neq (1,1)$, we have that $\mu(A_{1,1}\setminus A'_{1,1}) =
0$ and it follows from the above equality that also
$\mu^{\rm{F}}_{[d]}\big(\prod_{i=1}^d\big(\bigcap_{j=1}^{k_i}A'_{i,j}\big)\big)
= 0$, so an appeal to property P for the reduced collection of
up-sets completes the proof. \qed

\textbf{Remark}\quad The first very simple case treated by the above
proof is the only step in the whole of the present section that is
essentially absent from Tao's arguments in Sections 6 and 7
of~\cite{Tao07}.  Nevertheless, it seems to be essential for the
correct organization of the present argument, since we need to allow
for which of our sets are lifted under which coordinate projections
in the hypothesis that
$\mu_{[d]}^{\rm{F}}\big(\prod_{i=1}^d\big(\bigcap_{j=1}^{k_i}A_{i,j}\big)\big)
= 0$. \fin

\begin{lem}[Lifting under finitary generation]\label{lem:lift-gen}
Suppose that all up-sets in the collection $(\I_{i,j})_{i,j}$ have
depth at least $k$ and that among those of depth $k$ there are $\ell
\geq 1$ that are non-principal.  Then if property P holds for all
similar collections having at most $\ell - 1$ non-principal up-sets
of depth $k$, then it also holds for this collection.
\end{lem}

\textbf{Proof}\quad Let $\I_{i_1,j_1}$, $\I_{i_2,j_2}$, \ldots,
$\I_{i_\ell,j_\ell}$ be the non-principal up-sets of depth $k$, and
now in addition let $e_1$, $e_2$, \ldots, $e_r$ be all the members
of $\I_{i_\ell,j_\ell}$ of size $k$ (so, of course, $r \leq
\binom{d}{k}$). Once again we will assume for simplicity that
$(i_\ell,j_\ell) = (1,1)$.  We break our work into two further
steps.

\textbf{Step 1}\quad First consider the case of a collection
$(A_{i,j})_{i,j}$ such that for the set $A_{1,1}$, we can actually
find \emph{finite} subalgebras of sets $\B_s\in \Phi_{\{e_s\}}$ for
$s = 1,2,\ldots,r$ such that $A_{i_\ell,j_\ell} \in \B_1\vee
\B_2\vee \cdots \vee \B_r\vee\Phi_{\I_{1,1}\cap\binom{[d]}{\geq
k+1}}$ (so $A_{1,1}$ lies in one of our non-principal up-sets of
depth $k$, but it fails to lie in an up-set of depth $k+1$ only `up
to' finitely many additional generating sets). Choose $M\geq
\max_{s\leq r}|\B_s|$, so that we can certainly express
\[A_{1,1} = \bigcup_{m=1}^{M^r} (B_{m,1}\cap B_{m,2}\cap\cdots\cap B_{m,r}\cap C_m)\]
with $B_{m,s}\in \B_s$ for each $s\leq r$ and $C_m \in
\Phi_{\I_{1,1}\cap\binom{[d]}{\geq k+1}}$.  Inserting this
expression into the equation
\[\mu^{\rm{F}}_{[d]}\Big(\prod_{i=1}^d\Big(\bigcap_{j=1}^{k_i}A_{i,j}\Big)\Big)
= 0\] now gives that each of the $M^r$ individual sets
\[\Big((B_{m,1}\cap B_{m,2}\cap\cdots\cap B_{m,r}\cap C_m)\cap \bigcap_{j=2}^{k_1}A_{1,j}\Big)\times\prod_{i=2}^d \Big(\bigcap_{j=1}^{k_i}A_{i,j}\Big)\]
is $\mu^{\rm{F}}_{[d]}$-negligible.

Now consider the family of up-sets comprising the original
$\I_{i,j}$ if $i=2,3,\ldots,d$ and the collection $\langle
e_1\rangle$, $\langle e_2\rangle$, \ldots, $\langle e_r\rangle$,
$\I_{1,2}$, $\I_{1,3}$, \ldots, $\I_{1,k_1}$ corresponding to $i =
1$. We have broken the depth-$k$ non-principal up-set $\I_{1,1}$
into the higher-depth up-set $\I_{1,1}\cap\binom{[d]}{\geq k+1}$ and
the principal up-sets $\langle e_s\rangle$, and so there are only
$\ell-1$ minimal-depth non-principal up-sets in this new family.  It
is clear that for each $m\leq M^r$ the above product set is
associated to this family of up-sets, and so an inductive appeal to
property P for this family tells us that also
\[\mu\Big((B_{m,1}\cap B_{m,2}\cap\cdots\cap B_{m,r}\cap C_m)\cap \bigcap_{j=2}^{k_1}A_{1,j}\cap \bigcap_{i=2}^d \bigcap_{j=1}^{k_i}A_{i,j}\Big) = 0\]
for every $m\leq M^r$.  Since the union of these sets is just
$\bigcap_{i=1}^d\bigcap_{j=1}^{k_i} A_{i,j}$, this gives the desired
negligibility in this case.

\textbf{Step 2}\quad Now we return to the general case, which will
follow by a suitable limiting argument applied to the conclusion of
Step 1.  Since any $\Phi_{\{e\}}$ is countably separated, for each
$e$ with $|e| = k$ we can find an increasing sequence of finite
subalgebras $\B_{e,1} \subseteq \B_{e,2} \subseteq\ldots$ that
generates $\Phi_{\{e\}}$ up to $\mu$-negligible sets. In terms of
these define approximating sub-$\s$-algebras
\[\Xi^{(n)}_{i,j} := \Phi_{\I_{i,j}\cap \binom{[d]}{\geq k+1}}\vee \bigvee_{e \in \I_{i,j}\cap \binom{[d]}{k}} \B_{e,n},\]
so for each $\I_{i,j}$ these form an increasing family of
$\s$-algebras that generates $\Phi_{\I_{i,j}}$ up to
$\mu$-negligible sets (indeed, if $\I_{i,j}$ does not contain any
sets of the minimal depth $k$ then we simply have $\Xi^{(n)}_{i,j} =
\Phi_{\I_{i,j}}$ for all $n$).

Observe that by Corollary~\ref{cor:Fberg-struct}, for each $n$ we
have that $\Phi^{\rm{F}}_{\I_{1,1}}$ and $\bigvee_{(i,j)\neq
(1,1)}\pi_i^{-1}(\Xi^{(n)}_{i,j})$ are relatively independent over
$\pi_1^{-1}(\Xi^{(n)}_{1,1})$.

Given now a family of sets $(A_{i,j})_{i,j}$ associated to
$(\I_{i,j})_{i,j}$, for each $(i,j)$ the conditional expectations
$\sfE_\mu(1_{A_{i,j}}\,|\,\Xi^{(n)}_{i,j})$ form an almost surely
uniformly bounded martingale converging to $1_{A_{i,j}}$ in
$L^2(\mu)$. Letting $B^{(n)}_{i,j}:=
\{\sfE_\mu(1_{A_{i,j}}\,|\,\Xi^{(n)}_{i,j})
> 1-\delta\}$ for some small $\delta > 0$ (to be specified
momentarily), it is clear that we also have $\mu(A_{i,j}\triangle
B_{i,j}^{(n)}) \to 0$ as $n\to\infty$.  Let also
\[F := \prod_{i=1}^d\Big(\bigcap_{j=1}^{k_i}B^{(n)}_{i,j}\Big).\] We
now compute using the above-mentioned relative independence that
\begin{eqnarray*}
\mu^{\rm{F}}_{[d]}(F\setminus \pi_i^{-1}(A_{i,j})) &=& \int_{X^d}
\Big(\prod_{(i',j')}1_{B^{(n)}_{i',j'}}\circ\pi_{i'}\Big) - 1_{A_{i,j}}\circ\pi_i\cdot\Big(\prod_{(i',j')}1_{B^{(n)}_{i',j'}}\circ\pi_{i'}\Big)\,\d\mu^{\rm{F}}_{[d]}\\
&=& \int_{X^d} (1_{B^{(n)}_{i,j}\setminus A_{i,j}}\circ\pi_i)\cdot
\Big(\prod_{(i',j')\neq
(i,j)}1_{B^{(n)}_{i',j'}}\circ\pi_{i'}\Big)\,\d\mu^{\rm{F}}_{[d]}\\
&=& \int_{X^d} (\sfE_\mu(1_{B^{(n)}_{i,j}\setminus
A_{i,j}}\,|\,\Xi^{(n)}_{i,j})\circ\pi_i)\cdot
\Big(\prod_{(i',j')\neq
(i,j)}1_{B^{(n)}_{i',j'}}\circ\pi_{i'}\Big)\,\d\mu^{\rm{F}}_{[d]}
\end{eqnarray*}
for each pair $(i,j)$.

However, from the definition of $B^{(n)}_{i,j}$ we must have
\[\sfE_\mu(1_{B^{(n)}_{i,j}\setminus A_{i,j}}\,|\,\Xi^{(n)}_{i,j}) \leq \delta 1_{B^{(n)}_{i,j}}\]
almost surely, and therefore the above integral inequality implies
that \[\mu^{\rm{F}}_{[d]}(F\setminus \pi_i^{-1}(A_{i,j}))\leq
\delta\int_{X^d} (1_{B^{(n)}_{i,j}}\circ\pi_i)\cdot
\Big(\prod_{(i',j')\neq
(i,j)}1_{B^{(n)}_{i',j'}}\circ\pi_{i'}\Big)\,\d\mu^{\rm{F}}_{[d]} =
\delta \mu^{\rm{F}}_{[d]} (F).\]

From this we can estimate as follows: \[\mu^{\rm{F}}_{[d]}(F) \leq
\mu^{\rm{F}}_{[d]}\Big(\prod_{i=1}^d\Big(\bigcap_{j=1}^{k_i}A_{i,j}\Big)\Big)
+ \sum_{(i,j)}\mu^{\rm{F}}_{[d]}(F\setminus \pi_i^{-1}(A_{i,j}))
\leq 0 + \Big(\sum_{i=1}^dk_i\Big)\delta \mu^{\rm{F}}_{[d]} (F),\]
and so provided we chose $\delta < \big(\sum_{i=1}^dk_i\big)^{-1}$
we must in fact have $\mu^{\rm{F}}_{[d]} (F) = 0$.

We have now obtained sets $(B^{(n)}_{i,j})_{i,j}$ that are
associated to the family $(\I_{i,j})_{i,j}$ and satisfy the property
of lying in finitely-generated extensions of the relevant factors
corresponding to the members of the $\I_{i,j}$ of minimal size, and
so we can apply the result of Step 1 to deduce that
$\mu\big(\bigcap_{i=1}^d\bigcap_{j=1}^{k_i}B^{(n)}_{i,j}\big) = 0$.
It follows that
\[\mu\Big(\bigcap_{i=1}^d\bigcap_{j=1}^{k_i}A_{i,j}\Big) \leq \sum_{i,j}\mu(A_{i,j}\setminus B^{(n)}_{i,j}) \to 0\quad\quad\hbox{as }n\to\infty,\]
as required. \qed

\textbf{Proof of Proposition~\ref{prop:infremoval}}\quad We first
take as our base case $k_i = 1$ and $\I_{i,1} = \{[d]\}$ for each
$i=1,2,\ldots,d$.  In this case we know that for any $A \in
\Phi_{[d]}$ the pre-images $\pi_i^{-1}(A)$ are all equal up to
negligible sets, and so given $A_1$, $A_2$, \ldots, $A_d \in
\Phi_{[d]}$ we have $0 = \mu^{\rm{F}}_{[d]}(A_1\times
A_2\times\cdots\times A_d) = \mu(A_1\cap A_2\cap\cdots\cap A_d)$.

The remainder of the proof now just requires putting the preceding
lemmas into order to form an induction with three layers: if our
collection has any non-principal up-sets of minimal depth, then
Lemma~\ref{lem:lift-gen} allows us to reduce their number at the
expense only of introducing new principal up-sets of the same depth;
and having removed all the non-principal minimal-depth up-sets,
Lemma~\ref{lem:lift-rel-ind} enables us to remove also the principal
ones until we are left only with up-sets of increased minimal depth.
This completes the proof. \qed

\parskip 0pt

\bibliographystyle{abbrv}
\bibliography{newmultiSzem}

\begin{thebibliography}{10}

\bibitem{Aus--nonconv}
T.~Austin.
\newblock On the norm convergence of nonconventional ergodic averages.
\newblock To appear, \emph{{E}rgodic {T}heory {D}ynam. {S}ystems}.

\bibitem{ConLes88.1}
J.-P. Conze and E.~Lesigne.
\newblock Sur un th\'eor\`eme ergodique pour des mesures diagonales.
\newblock In {\em Probabilit\'es}, volume 1987 of {\em Publ. Inst. Rech. Math.
  Rennes}, pages 1--31. Univ. Rennes I, Rennes, 1988.

\bibitem{ConLes88.2}
J.-P. Conze and E.~Lesigne.
\newblock Sur un th\'eor\`eme ergodique pour des mesures diagonales.
\newblock {\em C. R. Acad. Sci. Paris S\'er. I Math.}, 306(12):491--493, 1988.

\bibitem{FraKra05}
N.~Frantzikinakis and B.~Kra.
\newblock Convergence of multiple ergodic averages for some commuting
  transformations.
\newblock {\em Ergodic Theory Dynam. Systems}, 25(3):799--809, 2005.

\bibitem{Fur77}
H.~Furstenberg.
\newblock Ergodic behaviour of diagonal measures and a theorem of {S}zemer\'edi
  on arithmetic progressions.
\newblock {\em J. d'Analyse Math.}, 31:204--256, 1977.

\bibitem{FurKat78}
H.~Furstenberg and Y.~Katznelson.
\newblock An ergodic {S}zemer\'edi {T}heorem for commuting transformations.
\newblock {\em J. d'Analyse Math.}, 34:275--291, 1978.

\bibitem{Gla03}
E.~Glasner.
\newblock {\em Ergodic {T}heory via {J}oinings}.
\newblock American {M}athematical {S}ociety, {P}rovidence, 2003.

\bibitem{Gow??}
W.~T. Gowers.
\newblock Hypergraph regularity and the multidimensional {S}zemer\'edi
  {T}heorem.
\newblock preprint.

\bibitem{Hos??}
B.~Host.
\newblock Ergodic seminorms for commuting transformations and applications.
\newblock Preprint.

\bibitem{HosKra01}
B.~Host and B.~Kra.
\newblock Convergence of {C}onze-{L}esigne averages.
\newblock {\em Ergodic Theory Dynam. Systems}, 21(2):493--509, 2001.

\bibitem{HosKra05}
B.~Host and B.~Kra.
\newblock Nonconventional ergodic averages and nilmanifolds.
\newblock {\em Ann. Math.}, 161(1):397--488, 2005.

\bibitem{NagRodSch07}
B.~Nagle, V.~R\"odl, and M.~Schacht.
\newblock The counting lemma for regular $k$-uniform hypergraphs.
\newblock {\em Random Structures and Algorithms}, to appear.

\bibitem{Tao07}
T.~Tao.
\newblock A correspondence principle between (hyper)graph theory and
  probability theory, and the (hyper)graph removal lemma.
\newblock {\em J. d'Analyse Math.}, 103:1--45, 2007.

\bibitem{Tao08(nonconv)}
T.~Tao.
\newblock Norm convergence of multiple ergodic averages for commuting
  transformations.
\newblock {\em Ergodic Theory Dynam. Systems}, 28:657--688, 2008.

\bibitem{Zha96}
Q.~Zhang.
\newblock On convergence of the averages {$(1/N)\sum\sp N\sb {n=1}f\sb 1(R\sp
  nx)f\sb 2(S\sp nx)f\sb 3(T\sp nx)$}.
\newblock {\em Monatsh. Math.}, 122(3):275--300, 1996.

\bibitem{Zie07}
T.~Ziegler.
\newblock Universal characteristic factors and {F}urstenberg averages.
\newblock {\em J. Amer. Math. Soc.}, 20(1):53--97 (electronic), 2007.

\end{thebibliography}

\vspace{10pt}

\small{\textsc{Department of Mathematics, University of California,
Los Angeles CA 90095-1555, USA}}

\vspace{5pt}

\small{Email: \verb|timaustin@math.ucla.edu|}

\vspace{5pt}

\small{URL: \verb|http://www.math.ucla.edu/~timaustin|}

\end{document}